\documentclass[12pt]{article}
\usepackage{amsthm,amsmath}
\usepackage{amscd}
 \usepackage{amssymb}
 \usepackage{hyperref}
 \usepackage[all]{xypic}
 \usepackage{mathtools}
 \usepackage[utf8]{inputenc}
 \usepackage{enumerate}
  \newtheorem{lemma}{Lemma}[section]
 \newtheorem{corollary}[lemma]{Corollary}
 
 \newtheorem{theorem}[lemma]{Theorem}

\usepackage[utf8]{inputenc} 

\def\includegraphics{}

\usepackage{tikz,xcolor,hyperref}

\definecolor{lime}{HTML}{A6CE39}
\DeclareRobustCommand{\orcidicon}{
	\begin{tikzpicture}
	\draw[lime, fill=lime] (0,0)
	circle[radius=0.16]
	node[white]{{\fontfamily{qag}\selectfont \tiny \.{I}D}};
	\end{tikzpicture}
	\hspace{-2mm}
}
\foreach \x in {A, ..., Z}{%
	\expandafter\xdef\csname orcid\x\endcsname{\noexpand\href{https://orcid.org/\csname orcidauthor\x\endcsname}{\noexpand\orcidicon}}
}

\usepackage{amscd}
\usepackage{amsmath,amssymb,mathrsfs}
\usepackage{hyperref}
\usepackage[all]{xypic}
\usepackage{times}
\usepackage{inputenc}
\usepackage{titlesec}

\allowdisplaybreaks[1]
\allowdisplaybreaks[4]

 \usepackage{amsthm}
 \usepackage{multirow}
\usepackage{array}
\usepackage{tabularx}
\usepackage[figuresright]{rotating}
\usepackage{bm}

\advance\headheight1.15pt

    \newcommand{\Aut}{\mathrm{Aut}}
  
      \newcommand{\calS}{\mathcal{S}}

    \newcommand{\GF}{\mathrm{GF}} \newcommand{\GL}{\mathrm{GL}}

  \newcommand{\PGL}{\mathrm{PGL}} \newcommand{\PGaL}{\mathrm{P\Gamma L}}     \newcommand{\PSL}{\mathrm{PSL}}

 \newcommand{\SL}{\mathrm{SL}}        \newcommand{\Sym}{\mathrm{Sym}}

\renewcommand{\mod}{\mathrm{mod}\,}

\newcommand{\Fix}{\mathrm{Fix}}

\begin{document}

\baselineskip=16pt
\title{\fontsize{15}{0}\selectfont  On the Cameron-Praeger conjecture for Steiner 6-designs}
\author{\fontsize{11}{0}\selectfont
	Ximin Wang$^{a}$,~~Zheng Huang$^{*,a}$,~~Weijun Liu$^{b,c}$
	\\
	\fontsize{10}{0}\it{} 
	\\ \fontsize{10}{0}\it{  $^{a}$School of Mathematics and Statistics, Hunan University of Science and Technology,}
	\\ \fontsize{10}{0}\it{ Xiangtan, Hunan 411201, China}
\\ \fontsize{10}{0}\it{  $^{b}$College of General Education, Guangdong University of Science and Technology,}
	\\ \fontsize{10}{0}\it{ Dongguan, Guangdong 523083, China}
\\ \fontsize{10}{0}\it{  $^{c}$School of Mathematics and Statistics, Central South University,}
	\\ \fontsize{10}{0}\it{ Changsha, Hunan 410083, China}
}

\date{}
\maketitle
\thispagestyle{empty}
\begin{center}

\textbf{\underline{ABSTRACT}}
\end{center}
\ \ \ \ \

A long-standing and far-reaching conjecture of Cameron and Praeger (1993) asserts that no non-trivial block-transitive $6$-$(v,k,\lambda)$ designs exist.
Huber (2010) reduced the Steiner case ($\lambda=1$) to a single infinite family of possible exceptions, where the automorphism group is $G=\PGaL(2,p^e)$ acting on the points of the projective line, with $p\in\{2,3\}$ and $e$ an odd prime power.
In this paper, we resolve this final case.
By combining orbit-stabiliser analysis for the action of specific involutions in $\PGaL(2,p^e)$ with the modular arithmetic constraints derived from the design parameters, we prove that no such Steiner 6-designs can exist for any admissible $q$.
This completes the proof that the Cameron-Praeger conjecture is true for the important case of non-trivial Steiner 6-designs.

\begin{flushleft}
\end{flushleft}
\begin{flushleft}
\end{flushleft}

\let\thefootnote\relax\footnote{$^*$Corresponding Author}
\let\thefootnote\relax\footnote{Ximin Wang: wxm240325@163.com}
\let\thefootnote\relax\footnote{ Zheng Huang: huangzhengmath@163.com}
\let\thefootnote\relax\footnote{ Weijun Liu: wjliu6210@126.com }
\vspace{1cm}

\section{Introduction} \label{sec1}

\ \ \  \ \ \ A $t$-$(v,k,\lambda)$ design  $\mathcal{D}= (\mathcal{P},\mathcal{B})$  is an incidence structure consisting of a set $\mathcal{P}$ of $v$ \emph{points} and a set $\mathcal{B}$ with $k$-subsets of $\mathcal{P}$ called \emph{blocks}  such that, each block in $\mathcal{B}$ has size $k$, and each $t$-subset of $\mathcal{P}$ lies in exactly $\lambda$ blocks from $\mathcal{B}$.
An {\it automorphism} of $\mathcal{D}$ is a permutation of the points which preserves the blocks.
The set of all automorphisms of $\mathcal{D}$ under the composition of permutations forms a group, denoted by $\Aut(\mathcal{D})$.
If $G\leq\mathrm{Aut}(\mathcal{D})$, then $G$ is called  an {\it automorphism group} of $\mathcal{D}$.
$\mathcal{D}$ is {\it $G$-block-transitive} if $G$ acts transitively on the blocks, and we also say $G$ is {\it  block-transitive}.

In 2002, renowned permutation group expert Praeger (a Fellow of the Australian Academy of Sciences) delivered a 45-minute invited lecture at the International Congress of Mathematicians, systematically introducing research on finite permutation groups and block designs.
To date, the classification of block designs with certain transitivity properties in their automorphism groups remains an active topic in algebraic combinatorics.


In 1993, Cameron and Praeger \cite{Cameron1} conjectured that there are no non-trivial block-transitive 6-designs.
The conjecture later became known as the {\bf Cameron-Praeger conjecture} and emerged as one of the major open problems in permutation group theory and combinatorial design theory.

In fact, the Cameron-Praeger conjecture is essentially a conjecture about the upper bound on $t$ for nontrivial block-transitive $t$-designs.
Since block-transitive 5-designs do exist, proving that no block-transitive 6-design exists would establish that the upper bound is 5.
In \cite{Cameron1}, Cameron and Praeger also gave an important result for the study of block-transitive designs.
They proved that for $t$-design $\mathcal{D}$ with $t\geq2$, if $G\leq\Aut(\mathcal{D})$ acts block-transitively on $\mathcal{D}$, then $G$ acts point $\left\lfloor\frac{t}{2}\right\rfloor$-homogeneously on $\mathcal{D}$.
According to this conclusion, when $t=6$, all $G$-block-transitive $6$-design are point $3$-homogeneously.
Therefore, research on the Cameron-Praeger conjecture will focus on finite 3-homogeneous permutation groups.

Up to now, progress on the Cameron-Praeger conjecture has been relatively slow.
Huber \cite{Huber1} treated  $6$-$(v,k,1)$ designs and reduced their possible block-transitive groups to particular projective semilinear groups in characteristics two and three.
Chen and Liu \cite{Chen} proved nonexistence when $k\leq100$, without restricting the index.
Tan, Liu and Chen \cite{Tan} extended the bound to $k\leq10000$ when the block-transitive group is almost simple.
More recently, Wang, Huang, Zhang and Liu \cite{Wang} proved that a block-transitive
$6$-$(v,k,2)$ design does not exist under the additional assumption $k\mid v$.
These results suggest that the case of small index provides a useful setting in which to examine the remaining structural obstructions.

Building upon Huber's work, this paper eliminates the exceptional cases in his theorem, thereby completing the proof of the nonexistence of nontrivial block-transitive Steiner 6-designs.
Our main theorems are as follows:

\begin{theorem}\label{Mainth}
There is no non-trivial block-transitive Steiner $6$-designs.
\end{theorem}

 \section{Preliminaries}\label{sec2}


\ \ \ \ \ \
Since this paper makes extensive use of design, finite fields and permutation groups, we first review some basic concepts related to these topics.

Let $\Omega$ be a set of $n$ points, and let $\Omega_{n}^{\{s\}}$ denote the family of all $s$-subsets of $\Omega$ for $1\leq s\leq n$.
Let $\mathcal{D}=(\mathcal{P},\mathcal{B})$ be a $t$-design.
The parameters of $\mathcal{D}$ is given by $(v,b,r,k)$, where $v$ is the number of points, $b$ is the number of blocks, $r$ is the number of blocks that incident with a given point and $k$ is the number of points of a block.
Let $\mathcal{P}^{\{s\}}$ denote the family of all $s$-subsets of $\mathcal{P}$.
For $1\leq s\leq t$, an \emph{$s$-flag} of $\mathcal{D}$ is a point-set block pair $\{S,B\}$, where $S\in \mathcal{P}^{\{s\}}$ and $S\subset B \in \mathcal{B}$.

Let $\GF(q)$ be a finite field of order $q$ and characteristic $p$, and $q=p^e$.
The projective line is denoted by $\GF(q)\cup\{\infty\}$.
$\GL(2,q)$ denotes the group of invertible two-by-two matrices over $\GF(q)$, and $\SL(2,q)$ is its subgroup of determinant-one matrices.
The quotients of $\GL(2,q)$ and $\SL(2,q)$ by their centers are denoted $\PGL(2,q)$ and  $\PSL(2,q)$, respectively.

Let $A=
\begin{pmatrix}
a & b \\
c & d
\end{pmatrix}$
be an invertible matrix over the field $\GF(q)$.
\begin{align*}
t_{a,b,c,d}: x\mapsto\frac{ax+b}{cx+d}
\end{align*}
is a fractional linear transformation of the projective line $\GF(q)\cup\{\infty\}$ that induced by $A$, and $\PGL(2,q)$ is the group of all fractional linear transformations.
For $\lambda\in\GF(q)$, $t_{a,b,c,d}=t_{\lambda a,\lambda b,\lambda c,\lambda d}$.
And $\PSL(2,q)=\{t_{a,b,c,d}\mid ad-bc~is~a~square\}$.
\begin{align*}
\Aut(\PSL(2,q))=\PGaL(2,q)=\PGL(2,q)\rtimes\langle \tau_{\alpha}\rangle,
\end{align*}
where $\tau_{\alpha}\in\Sym(\GF(q)\cup\{\infty\})\cong \calS_{q+1}$ with order $e$ is a permutation induced by the Frobenius automorphism $\alpha: \GF(p^e)\rightarrow\GF(p^e), x\mapsto x^p$.


\begin{lemma}{\rm\cite{Guralnick}}\label{lem2.1}
$\PSL(2,q)$ contains a unique conjugacy class of involutions.
For $q$ even, the centralizer of an involution of $\PGL(2,q)$ is a Sylow 2-subgroup.
For $q$ odd, the centralizer in $\PGL(2,q)$ of an involution of $\PGL(2,q)$ is dihedral of order either $2(q+1)$ or $2(q-1)$, and the centralizer in $\PSL(2,q)$ of this involution is dihedral of half the size.
\end{lemma}

\begin{lemma}{\rm\cite{Burness}}\label{lem2.2}
Let $G$ be a permutation group on a finite set $\Omega$ and $\alpha\in\Omega$.
For any $g\in G$, we have
\begin{align*}
\mid \Fix_{\Omega}(g)\mid=\frac{\mid C_{G}(g)\mid\mid g^{G}\cap G_{\alpha}\mid}{\mid G_{\alpha}\mid}\leq\mid C_{G}(g)\mid
\end{align*}
\end{lemma}

The following lemma  can be derived from equation (4.2) in Theorem 4.1 of \cite{Diaconis} and will facilitate the calculation of the design parameters.

\begin{lemma}{\rm\cite{Diaconis}}\label{lem2.3}
Let $z\in \calS_n$ with $z^2=1$, and let $f$ be the number of fixed points of $z$ in its usual action on $n$ points.
If $f\leq2$, then the number of fixed points of $z$ acting on $\Omega_{n}^{2s}$ is ${\left\lfloor n/2\right\rfloor \choose s}$.
That is, the number of invariant 2-subsets under the natural action of $z$ on the points is
$\Omega_{n}^{2s}$ is ${\left\lfloor n/2\right\rfloor \choose s}$.
\end{lemma}

Next, we give some well known fact of designs.

\begin{lemma}{\rm\cite{Beth}}\label{lem2.4}
If $\mathcal{D}=(\mathcal{P},\mathcal{B})$ is a non-trivial $t$-$(v,k,\lambda)$ design with $t\geq2$, then
\begin{enumerate}[{\rm(1)}]
\item $bk=vr$;

\item ${v \choose t }\lambda=b{k \choose t}$;

\item for a positive integer $1\leq s\leq t$, a $\mathcal{D}$ also is an $s$-$(v,k,\lambda_s)$ design, where $\lambda_s=\lambda\frac{\binom{v-s}{t-s}}{\binom{k-s}{t-s}}$.
\end{enumerate}
\end{lemma}

\begin{lemma}{\rm\cite{Cameron1}}\label{lem2.5}
Let $\mathcal{D}$ be a $t$-design with $t\geq2$.
Then the following holds:
\begin{enumerate}[{\rm(i)}]

\item If $G\leq\Aut(\mathcal{D})$ acts block-transitively on $\mathcal{D}$,
then $G$ acts point $\left\lfloor\frac{t}{2}\right\rfloor$-homogeneously on $\mathcal{D}$.

\item If $G\leq\Aut(\mathcal{D})$ acts flag-transitively on $\mathcal{D}$,
then $G$ acts point $\left\lfloor\frac{t+1}{2}\right\rfloor$-homogeneously on $\mathcal{D}$.
\end{enumerate}
\end{lemma}

\begin{lemma}\label{lem3.2}
Let $\mathcal{D}=(\mathcal{P},\mathcal{B})$ be a non-trivial Steiner $t$-design.
Then the following holds:
\begin{enumerate}[{\rm(i)}]

\item {\rm(Tits, 1964 \cite{Tits})} $v\geq(t+1)(k-t+1)$;

\item {\rm(Cameron, 1976 \cite{Cameron0})} $v-t+1\geq(k-t+2)(k-t+1)$ for $t>2$.
If equality holds, then $(t,k,v)=(3,4,8),~(3,6,22),~(3,12,112),~(4,7,23)$ or $(5,8,24)$.
\end{enumerate}
\end{lemma}
From (ii), Huber deduce the following upper bound for the positive integer $k$.
\begin{lemma}{\rm\cite{Huber1}}\label{lem3.3}
Let $\mathcal{D}=(\mathcal{P},\mathcal{B})$ be a non-trivial Steiner $6$-design.
Then
\begin{align*}
k\leq\left\lfloor \sqrt{v-\frac{19}{4}}+\frac{9 }{2}\right\rfloor.
\end{align*}
\end{lemma}

Furthermore, according to Lemma \ref{lem3.3}, for certain specific numbers of points, we have the following strict upper bound for the positive integer $k$.
\begin{corollary}\label{cor3.4}
Let $\mathcal{D}=(\mathcal{P},\mathcal{B})$ be a non-trivial Steiner $6$-design.
If $v=q+1$, we have
\begin{align*}
k\leq\left\lfloor \sqrt{q-\frac{15}{4}}+\frac{9 }{2}\right\rfloor<\left\lfloor \sqrt{q}+\frac{9 }{2}\right\rfloor.
\end{align*}
\end{corollary}

\begin{lemma}\label{lem3.5}
Let $\mathcal{D}=(\mathcal{P},\mathcal{B})$ be a non-trivial Steiner $2s$-design.
If $G\leq\Aut(\mathcal{D})$ acts block-transitively on $\mathcal{D}$, and for a involution $g\in G$, it fixes at most two points in the point set $\mathcal{P}$.
Then
\begin{align*}
|\Fix_{\mathcal{B}}(g)|=\frac{{\left\lfloor v/2\right\rfloor \choose s}}{{\left\lfloor k/2\right\rfloor \choose s}}\leq |C_G(g)|.
\end{align*}
\end{lemma}

\begin{proof}
Suppose that $\mathcal{D}=(\mathcal{P},\mathcal{B})$ is a $G$-block-transitive design.
Let $g\in G$ be an involution that fixes at most two points of $\mathcal{P}$.
According to Lemma \ref{lem2.3}, the number of invariant $2s$-subsets under $g$ is ${\left\lfloor v/2\right\rfloor \choose s}$.
Let $S$ be be one of these invariant $2s$-subset.
By the definition of Steiner design, there exist a unique block $B$ such that $S\in B$.
Then $(S,B)$ is a $2s$-flag of $\mathcal{D}$.
Since $g\in G$, we have $(S,B)^g=(S^g,B^g)=(S,B^g)$, where $(S,B^g)$ also a $2s$-flag of $\mathcal{D}$.
Hence $B=B^g$.
On the other hand, for some $B\in\mathcal{B}$, if $B^g=B$, then $g|_{B}$ also is a involution and $g|_{B}$ acts on $B$ at most two fixed points, where $g|_{B}$  is the permutation of $B$ induced by restricting $g$ to $B$.
Thus, by Lemma \ref{lem2.3}, the number of invariant $2s$-subset that $g|_{B}$ acts on $B$ is ${\left\lfloor k/2\right\rfloor \choose s}$.

Now, we count the number of $2s$-flag that fixed by $g$ and find that
\begin{align*}
{\left\lfloor v/2\right\rfloor \choose s}=|\Fix_{\mathcal{B}}(g)|{\left\lfloor k/2\right\rfloor \choose s}.
\end{align*}
By the definition of non-trivial design, we have $k>2s$.
So
\begin{align*}
{\left\lfloor k/2\right\rfloor \choose s}>0.
\end{align*}
Therefore, together with Lemma \ref{lem2.2}, we obtain
\begin{align*}
|\Fix_{\mathcal{B}}(g)|=\frac{{\left\lfloor v/2\right\rfloor \choose s}}{{\left\lfloor k/2\right\rfloor \choose s}}\leq |C_G(g)|.
\end{align*}
\end{proof}

\begin{theorem}{\rm\cite{Huber1}}\label{th3.6}
Let $\mathcal{D}=(\mathcal{P},\mathcal{B})$ be a non-trivial Steiner $6$-design.
Then $G\leq\Aut(\mathcal{D})$ cannot act block-transitively on $\mathcal{D}$, except possibly when $G=\PGaL(2,p^e)$ with $p=2$	or 3 and $e$ is an odd prime power.
\end{theorem}

\section{The proof of Main Theorem \ref{Mainth}}\label{sec3}

\ \ \ \ \ \  In this section, we complete the proof of the non-existence of nontrivial block-transitive Steiner 6-designs.

Let $\mathcal{D}=(\mathcal{P},\mathcal{B})$ be a non-trivial Steiner $6$-design.
According to Huber's result (Theorem \ref{th3.6}), it remains to consider only the two exceptional cases.
Therefore, in proving Theorem \ref{Mainth} below, we may assume without loss of generality that $\mathcal{P}=\GF(q)\cup\{\infty\}$, and then $v=q+1$.
$G=\PGaL(2,2^e)$  or  $G=\PGaL(2,3^e)$ with $e$ is an odd prime power.

First, for the above two types of groups, by studying the elements of the groups, we obtain the following conclusions:

\begin{lemma}\label{lem3.7}
\begin{enumerate}[{\rm(1)}]

\item If $q=2^e$, let $\omega(x)=x+1$.
Then, $\omega\in \PSL(2,q)$ with $|\Fix_{\mathcal{P}}(\omega)|=1$, and
\begin{align*}
|C_{\PGaL(2,q)}(\omega)|=qe.
\end{align*}

\item If $q=3^e$ with $e$ is an odd prime power, let $\omega(x)=-1/x$.
Then, $\omega\in \PSL(2,q)$ with $|\Fix_{\mathcal{P}}(\omega)|=0$, and
\begin{align*}
|C_{\PGaL(2,q)}(\omega)|=2(q+1)e.
\end{align*}
\end{enumerate}
\end{lemma}

\begin{proof}
For (1):  Let $q=2^e$ and $A=
\begin{pmatrix}
1 & 1 \\
0 & 1
\end{pmatrix}$.
$A$ is an invertible matrix over the field $\GF(2^e)$.
Then
$$\omega=t_{1,1,0,1}\in\PSL(2,q).$$
Note that $A\neq I$ and $A^2=I$, so $\omega$ is an involution.
For any $x\in\GF(2^e)\cup\{\infty\}$, if $x\neq\infty$, then $x+1=x$ is impossible.
Hence
$$|\Fix_{\mathcal{P}}(\omega)|=1.$$
By Lemma \ref{lem2.1}, we have $|C_{\PGL(2,q)}(\omega)|=q$.
As $\PGaL(2,q)=\PGL(2,q)\rtimes\langle \tau_{\alpha}\rangle$, for any $\tau_{\alpha}^i\in\langle \tau_{\alpha}\rangle$, by finite field theory, it follows that $\omega\tau_{\alpha}^i=\tau_{\alpha}^i\omega$.
Then, for any $g\in\PGaL(2,q)$, $\omega g=g\omega$ if and only if $g=u\tau_{\alpha}^i$ with $\omega u=u\omega$.
So, if $u\in\PGL(2,q)$, we have $u\in C_{\PGL(2,q)}$.
Therefore,
\begin{align*}
|C_{\PGaL(2,q)}(\omega)|=|C_{\PGL(2,q)}(\omega)||\langle \tau_{\alpha}\rangle|=qe.
\end{align*}

For (2):  Let $q=3^e$ and $A=
\begin{pmatrix}
0 & -1 \\
1 & 0
\end{pmatrix}$.
$A$ is an invertible matrix over the field $\GF(3^e)$ with $e$ is an odd prime power.
Then
$$\omega=t_{0,-1,1,0}\in\PSL(2,q).$$
Note that $A\neq I$ and $A^2=-I$, so $\omega$ is an non-identity involution.
On the one hand, we have $\omega(0)=\infty$ and $\omega(\infty)=0$.
On the other hand, for $x\in\GF(3^e)\cup\{\infty\}$, if $\omega(x)=x$, we have $x^2=-1$.
Then $x$ is an element of $GF(q)^*$ with order 4.
According to the Lagrange's Theorem of finite group, we have $4\mid q-1$, which contradicts $q=3^e\equiv3(\mod4)$.
Thus
$$|\Fix_{\mathcal{P}}(\omega)|=0.$$
By Lemma \ref{lem2.1}, we have $|C_{\PGL(2,q)}(\omega)|=2(q+1)$.
A similar analysis to (1) gives that for any $g=u\tau_{\alpha}^i=\in\PGaL(2,q)$ with $u\in\PGL(2,q)$, $\omega g=g\omega$ if and only if $\omega u=u\omega$.
So $u\in C_{\PGL(2,q)}$.
Therefore,
\begin{align*}
|C_{\PGaL(2,q)}(\omega)|=|C_{\PGL(2,q)}(\omega)||\langle \tau_{\alpha}\rangle|=2(q+1)e.
\end{align*}
\end{proof}

\begin{lemma}\label{lem3.8}
Let $\mathcal{D}=(\mathcal{P},\mathcal{B})$ be a non-trivial $\PGaL(2,p^e)$-block-transitive $6$-$(q+1,k,1)$ design.
Then
\begin{enumerate}[{\rm(1)}]

\item if $q=2^e$, we have
\begin{align*}
(q-2)(q-4)\leq 8e(\left\lfloor k/2\right\rfloor)(\left\lfloor k/2\right\rfloor-1)(\left\lfloor k/2\right\rfloor-2).
\end{align*}

\item if $q=3^e$ with $e$ is an odd prime power, we have
\begin{align*}
(q-1)(q-3)\leq 16e(\left\lfloor k/2\right\rfloor)(\left\lfloor k/2\right\rfloor-1)(\left\lfloor k/2\right\rfloor-2).
\end{align*}
\end{enumerate}
\end{lemma}
\begin{proof}
According to Lemma \ref{lem3.7}, for both case, there exist an involution $\omega\in\PGaL(2,p^e)$ with
$|\Fix_{\mathcal{P}}(\omega)|\leq2$.
Then, by Lemma \ref{lem3.5}, we have
\begin{align}\label{eq1}
|\Fix_{\mathcal{B}}(\omega)|=\frac{{\left\lfloor v/2\right\rfloor \choose 3}}{{\left\lfloor k/2\right\rfloor \choose 3}}\leq |C_{\PGaL(2,p^e)}(\omega)|.
\end{align}
If $q=2^e$, then $v=q+1=2^e+1$, $\left\lfloor v/2\right\rfloor=\frac{q}{2}$ and $|C_{\PGaL(2,p^e)}(\omega)|=qe$.
Hence, according to equation (\ref{eq1}), we have
\begin{align*}
\frac{q}{2}(\frac{q}{2}-1)(\frac{q}{2}-2)\leq qe\left\lfloor k/2\right\rfloor(\left\lfloor k/2\right\rfloor-1)(\left\lfloor k/2\right\rfloor-2).
\end{align*}
Thus,
\begin{align*}
(q-2)(q-4)\leq 8e(\left\lfloor k/2\right\rfloor)(\left\lfloor k/2\right\rfloor-1)(\left\lfloor k/2\right\rfloor-2).
\end{align*}

If $q=3^e$ with $e$ is an odd prime power, then $v=q+1=3^e+1$ is even,
$\left\lfloor v/2\right\rfloor=\frac{q+1}{2}$ and $|C_{\PGaL(2,p^e)}(\omega)|=2(q+1)e$.
According to equation (\ref{eq1}), we have
\begin{align*}
\frac{q+1}{2}(\frac{q+1}{2}-1)(\frac{q+1}{2}-2)\leq 2(q+1)e\left\lfloor k/2\right\rfloor(\left\lfloor k/2\right\rfloor-1)(\left\lfloor k/2\right\rfloor-2).
\end{align*}
Thus,
\begin{align*}
(q-1)(q-3)\leq 16e(\left\lfloor k/2\right\rfloor)(\left\lfloor k/2\right\rfloor-1)(\left\lfloor k/2\right\rfloor-2).
\end{align*}

\end{proof}

For convenience in subsequent writing and computation, we denote $h=2\left\lfloor k/2\right\rfloor$ and $f(h)=h(h-2)(h-4)$.
Thus, the following corollary can be derived from Lemma \ref{lem3.8}.

\begin{corollary}\label{cor3.9}
Let $\mathcal{D}=(\mathcal{P},\mathcal{B})$ be a non-trivial $\PGaL(2,p^e)$-block-transitive $6$-$(q+1,k,1)$ design.
Then
\begin{enumerate}[{\rm(1)}]

\item if $q=2^e$, we have
\begin{align}\label{eq2}
(q-2)(q-4)\leq ef(h).
\end{align}

\item if $q=3^e$ with $e$ is an odd prime power, we have
\begin{align}\label{eq3}
(q-1)(q-3)\leq 2ef(h).
\end{align}
\end{enumerate}
Moreover, for both cases, we have $(q-2)(q-4)\leq 2ek^3$.
\end{corollary}

\begin{proof}
As $h=2\left\lfloor k/2\right\rfloor$, we obtain
\begin{align*}
8(\left\lfloor k/2\right\rfloor)(\left\lfloor k/2\right\rfloor-1)(\left\lfloor k/2\right\rfloor-2)=h(h-2)(h-4)=f(h)
\end{align*}
Therefore, inequalities (2) and (3) follow directly from Lemma \ref{lem3.8}.
Note that $h\leq k$, it follows that $f(h)\leq k^3$.
Furthermore, when $q>8$, the inequality $(q-1)(q-3)\geq(q-2)(q-4)$ always holds.
Thus, for both cases, we have $(q-2)(q-4)\leq 2ek^3$.
\end{proof}

Next, we obtain a constrained upper bound for e through the parameter constraints derived from the properties of involutions in the group.
\begin{lemma}\label{lem3.10}
Let $\mathcal{D}=(\mathcal{P},\mathcal{B})$ be a non-trivial Steiner $6$-design.
For $G=\PGaL(2,p^e)$ with $p=2$	or 3 and $e$ is an odd prime power, if $\mathcal{D}$ is $G$-block-transitive, then $e<13$.
\end{lemma}

\begin{proof}
Let $\mathcal{D}=(\mathcal{P},\mathcal{B})$ be a non-trivial $\PGaL(2,p^e)$-block-transitive Steiner $6$-design with $p=2$	or 3 and $e$ is an odd prime power.
First, we assume that $e\geq13$, then $q\geq2^{13}$ and $\sqrt{q}>18$.
On the one hand, according to Corollary \ref{cor3.4}, we have $k< \sqrt{q}+\frac{9 }{2}$.
Since $\sqrt{q}>18$, we obtain
\begin{align*}
k\leq \frac{5}{4}\sqrt{q}.
\end{align*}
Then, by Corollary \ref{cor3.9}, we always have
\begin{align*}
(q-2)(q-4)\leq 2ek^3< 2e(\frac{5}{4}\sqrt{q})^3.
\end{align*}
On the other hand, as $q\geq2^{13}$, it follows that
\begin{align*}
(q-2)(q-4)\geq \frac{3}{4}q^2.
\end{align*}
Hence $\frac{3}{4}q^2<2e(\frac{5}{4}\sqrt{q})^3$, which means
$$\sqrt{q}<\frac{125}{32}e.$$
Since $q\geq2^{13}$, we have
\begin{align}\label{eq4}
2^{\frac{e}{2}}<\frac{125}{32}e.
\end{align}
Let $g(x)=2^{\frac{x}{2}}-\frac{125}{32}x$.
For $x\geq13$, we have $g'(x)>0$.
So for $x\geq13$, $g(x)\geq g(13)>0$.
Thus, for $e\geq13$, inequality (\ref{eq4}) cannot hold, which yields a contradiction.
Therefore, $e<13$.
\end{proof}

\begin{lemma}\label{lem3.11}
Let $\mathcal{D}=(\mathcal{P},\mathcal{B})$ be a non-trivial Steiner $6$-design.
Assume $G=\PGaL(2,p^e)$ with $p=2$	or 3 and $e$ is an odd prime power.
Then, there is no $G$-block-transitive designs $\mathcal{D}$ with $e<13$.
\end{lemma}

\begin{proof}
As $e$ is an odd prime power, we obtain $e=3,5,7,9,11$.
And by non-triviality, we have $k\geq7$.
According to Corollary \ref{cor3.4}, we have
\begin{align*}
k\leq \sqrt{q}+\frac{9}{2}.
\end{align*}
Hence, for each given $e$, $q$ can only take the value $2^e$ or $3^e$, and the corresponding values of $k$ in each case can be enumerated one by one.
According to Corollary \ref{cor3.9}, for a fixed value of $q$, we have $(q-1)(q-3)\leq ef(h)$ or $(q-2)(q-4)\leq2ef(h)$.
This gives a constraint between $q$ and $k$.
Then, we list each case one by one below.
If $q=2$, we have

\begin{center}
\begin{tabular}{c|c|c|c}
\hline
$(e,q)$ &The range of $k$& $(q-2)(q-4)$& $\max(ef(h))$\\
\hline
(3,8) & 7 & 24 &144\\
(5,32) & 7--9 & 840&960\\
(7,128) &7--15 &15624 &11760\\
(9,512) & 7--27 &259080 &123552\\
(11,2048) & 7--49 & 4182024&1068872\\
\hline
\end{tabular}
\end{center}
If $q=3$, we have

\begin{center}
\begin{tabular}{c|c|c|c}
\hline
$(e,q)$ &The range of $k$& $(q-2)(q-4)$& $\max(ef(h))$\\
\hline
(3,27) & 7--9 & 624 &1152\\
(5,243) & 7--19 & 58080&40320\\
(7,2187) &7--51 &4774224 &1545600\\
(9,19683) & 7--144 &387341760 &51528960\\
(11,177147) & 7--425 & 31380351024&1653294720\\
\hline
\end{tabular}
\end{center}
From the two tables above, it can be seen that only $(e,q)=(3,8),(5,32),(3,27)$ satisfies $(q-1)(q-3)\leq \max(ef(h))$ or $(q-2)(q-4)\leq\max(2ef(h))$.
Hence, all other cases can be ruled out by the parameter constraints given in Corollary \ref{cor3.9}.

Note that, by Lemma \ref{lem2.4} (3), for $s=0,1,\cdots, 5$, the parameter $\lambda_s=\lambda\frac{\binom{v-s}{t-s}}{\binom{k-s}{t-s}}$  is an integer.
Now, we discuss the remaining three cases one by one around this property.
For $(e,q)=(3,8)$, then $v=q+1=9$ and $k=7$, we obtain that $\lambda_1=\frac{28}{3}$ is not an integer.
For $(e,q)=(5,32)$, then $v=q+1=33$ and $k=7, 8$ or $9$.
If $k=7$, we have we $\lambda_4=\frac{406}{3}$;
If $k=7$, we have we $\lambda_5=\frac{28}{3}$;
If $k=7$, we have we $\lambda_4=\frac{406}{10}$.
For $(e,q)=(3,27)$, then $v=q+1=28$ and $k=7, 8$ or $9$.
If $k=7$, we have we $\lambda_5=\frac{23}{2}$;
If $k=7$, we have we $\lambda_5=\frac{23}{3}$;
If $k=7$, we have we $\lambda_5=\frac{23}{4}$.
Hence, all three remaining cases can be eliminated.
Therefore, there is no non-trivial $\PGaL(2,p^e)$-block-transitive Steiner $6$-design with $e<13$.
\end{proof}

In summary, combining Lemma \ref{lem3.10} and Lemma \ref{lem3.11} together with Theorem \ref{th3.6} completes the proof of Theorem \ref{Mainth}.
This also shows that the Cameron-Praeger conjecture holds for the case of Steiner designs.


\section*{Conflicts of Interest}
The authors declare that there is no conflict of interest.

\section*{Funding}
Z. Huang is supported by the National Natural Science Foundation of China (Grant No. 12401020) and the Natural Science Foundation of Hunan Province (Grant No.2026JJ60324). 
W.J. Liu is supported by the National Natural Science Foundation of China(Grant No. 12471022).

\section*{Availability of data and materials}
Data sharing not applicable to this article as no datasets were generated or analysed during
the current study: the article describes entirely theoretical research.




\pdfbookmark[1]{ References}{4}
 \end{document}